\documentclass[12pt]{article}
\usepackage{MnSymbol}
\usepackage{amsmath,amsfonts,
verbatim}

\newcommand{\A}{\mathbb{A}}
\newcommand{\Q}{{\mathbb Q}}
\newcommand{\C}{{\mathbb C}}
\newcommand{\N}{{\mathbb N}}

\newcommand{\R}{{\mathbb R}}

\newcommand{\F}{\mathbb{F}}
\newcommand{\U}{{\mathbb U}}
\newcommand{\X}{{\mathbb X}}

\newcommand{\pr}{{\rm pr}}

\newcommand{\la}{\langle}
\newcommand{\ra}{\rangle}
\newtheorem{pkt}{}[section]  
\newcommand{\bpk}{\begin{pkt}\rm }  
\newcommand{\epk}{\end{pkt}} 
\newcommand{\inv}{^{-1}}   

\newcommand{\be}{\begin{equation}}  
\newcommand{\ee}{\end{equation}}

\newcommand{\dcl}{\mathrm{dcl}}

\newcommand{\LL}{\mathcal{L}}
\newcommand{\kk}{\mathrm{k}}
\newcommand{\G}{\mathrm{G}}
\newcommand{\E}{\mathrm{E}}

\usepackage{graphicx}

\usepackage{tikz}

\title{A topological $L_{\omega_1,\omega}$-invariant.}
\author{B.Zilber}

\begin{document}
\maketitle
\abstract{We suggest to look at formal sentences describing complex algebraic varieties together with their universal covers as topological invariants. We prove that for abelian varieties and Shimura varieties this is indeed a complete invariant, i.e. it determines the variety up to complex conjugation. }

\section{Introduction}
The  universal (and more general) cover of a  complex algebraic variety as an abstract structure has been a subject of study in model theory in recent years, see e.g. \cite{MBays}, \cite{BHP}, \cite{DH}, \cite{Eter}. More specifically, given a smooth complex variety $\X$ and its cover $\pi: \U\to \X$ presented in a natural universal language,
the aim has been to find  an $L_{\omega_1,\omega}$-sentence $\Sigma_\X$ which axiomatises  $\pi: \U\to \X$  categorically, that is  determines the isomorphism type of an uncountable model of $\Sigma_\X$ uniquely in a given cardinality, up to an abstract isomorphism (ignoring topology).
This proved to be the case in the classes where one has enough knowledge  of the Galois action and arithmetic of $\X,$ see references above.

We suggest here that such a  $\Sigma_\X$ can serve as a complete topological invariant of respective $\X.$ More precisely, assuming that in a given class of varieties $\X$ such  a categorical $L_{\omega_1,\omega}$-sentence $\Sigma_\X$ exists, we expect:

{\em If, for a complex variety $\X',$  its cover structure  $\pi': \U'\to \X'$ is a model of the same sentence $\Sigma_\X$ then there is a biholomorphic isomorphism either between $\X$ and $\X'$ or between $\X^c,$ the complex conjugate of $\X$, and $\X'.$ }

(Compare this with speculations in a similar vein  presented by M.Gavrilovich in \cite{Misha}.) 

In this note  the cases of abelian varieties and Shimura varieties are elaborated.  

I am grateful to Chris Daw and Martin Bays for valuable comments and corrections.

\section{The language of covers}

\bpk \label{data} Let $\kk\subseteq \C,$ a countable subfield,  $\{ \X_i: i\in I\}$ a  collection of
non-singular irreducible complex algebraic varieties (of $\dim >0$) defined over $\kk$ and  $I:=(I,\ge )$ a lattice with the minimal element $\emptyset$ determined by unramified $\kk$-rational epimorphisms $\pr_{i',i}:\X_{i'}\to \X_i,$ for $i'\ge i.$ Let  $\U(\C)$ be
a connected complex manifold and $\{ \pi_i: i\in I\}$ 
 a collection of holomorphic covering maps  (local bi-holomorphismss)     
$$\pi_i: \U(\C)\twoheadrightarrow \X_i(\C),\ \ \pr_{i',i}\circ \pi_{i'}=\pi_i.$$
as illustrated by the picture:
\begin{center}

 \begin{tikzpicture}[scale=1]

\draw (12.2,4) node {$\X_n(\C)$};

 \draw (10.9,4.9) node {$\pr_{m,n}$};

\draw (9.6,5.5) node {$\X_m(\C)$};

\draw[->][line width=0.3mm] (9.5, 5.3)--(11.8, 4);

\draw[->][line width=0.3mm](9.5,4.8)--(9.9,3.1);
\draw[->][line width=0.3mm](12.1,3.9)--(10.1,3.1);

\draw (10,2.9) node {$\X(\C)$};
\draw (11.5,3.4) node {$\pr_{n,0}$};
\draw (9.3,4) node {$\pr_{m,0}$};

\draw (10,8.5) node {$\mathbb{U}(\C)$};
\draw[->][line width=0.3mm] (9.8,8 )--(9.5, 6);
\draw[->][line width=0.3mm] (10.2,8 )--(12.1, 4.3);
\draw (9.4,7) node {$\pi_m$};
\draw (11.,6) node {$\pi_n$};
\draw[->][line width=0.2mm] (10,8 )--(10, 3.1);



 
\end{tikzpicture}

\end{center}

\epk
\bpk \label{defS}
The langauge of covers   which we use is introduced  in full generality
in \cite{OAEC} based on the simple idea first formulated by M.Gavrilovich in \cite{Misha}.   

\medskip

{\bf Definition.} Given $S\subset \U^m$  
we say that $S$ is {\bf  $\LL$-primitive}
 if   there are $I_S\subseteq I$ and Zariski closed $Z_i\subseteq \X_i^m,$ $i\in I_S,$  defined over $\kk,$ 
such that $$S= \bigcap_{i\in I_S} \pi_i\inv(Z_i).$$ 
 
\medskip

It is easy to see that the one-sorted structure $(\U; \LL),$ interprets any finite cover $\X_i$ of $\X$ as the structure with $m$-ary relations given by Zariski closed subsets of $\X_i^m$ over $\kk.$

\epk

\bpk
{\bf Examples.}

1. Let $x_0\in \X_i(\kk)$ be a rational point. Then $\pi_i\inv(x_0)$ is definable over $\kk.$ Each point $u\in \pi_i\inv(x_0)$ is type-definable over $\bar{\kk},$ the algebraic closure of $\kk.$

2. Let $\Gamma$ be the fundamental group of $\X(\C)$ acting on $\U,$
$\gamma\in \Gamma.$ Then 
 $$\mathrm{graph}(\gamma)=\{ \la u,\gamma\cdot u\ra: u\in \U\}$$
is an analytic  component of $\pi\inv(\mathrm{Diag}(\X\times \X)),$ the pre-image of the diagonal. It is type-definable in the language $\LL$  using parameters from algebraic closure of $\kk.$ In fact, it is $\L_{\omega_1,\omega}$-definable using parameters from $\U.$


3. Let $\U(\C)=\C$ and $\X$ an elliptic curve over $\kk$ or $\mathbb{G}_m,$ $\pi=\exp,$ the respective exponentiation. By definition $\exp(0)$
is a $\kk$-rational point.  It follows (using the fact that the $\pi(u)$ can be taken just to be $\exp(\frac{u}{n})$)
 that $\{ 0\}$ is  definable. And so is the 
the ternary relation $u_1+u_2-u_3=0$ on $\C.$

4. For the case of elliptic curves, the group $\Lambda\otimes \Q$ is endowed with
an alternating bilinear form corresponding to the {\em Weil pairing} on the torsion of the elliptic curve. It induces 
  a definable relation on $\Lambda$ which is described in \cite{MBays}, 4.4.

\epk
\bpk 
   Structures in this language were studied in \cite{OAEC} under the assumptions of o-minimal representation of the cover (assumptions 2.2 in \cite{OAEC}). These are obviously satisfied  when $\X(\C)$ is compact.

\epk
\section{Abelian varities}\label{s1}
Recall that by \cite{MBays} for each elliptic non-CM curve $\E$  over a field $\kk\subset \C$ there is a 
 categorical $L_{\omega_1,\omega}$-sentence $\Sigma_\E$ in the language with parameters $C$ interpreted as elements of $\kk,$ which models the cover structure $(\pi:\C \to \E(\C)):$
  $$(\pi:\C \to\E(\C))\vDash \Sigma_\E.$$

The same is true for the general case of abelian varieties $\A$ by \cite{BHP}. The statement $\Sigma_\A(\bar{\lambda})$ there is proved over the language naming generators $\bar{\lambda}$ of $\Lambda$ but this can 
be eliminated by taking instead
$$\Sigma^\sharp_\A:= \exists \bar{\lambda}\Sigma_\A(\bar{\lambda})\ \&
\ \mathrm{qftp}(\pi(\Q\otimes \Lambda)).$$
(The $\Sigma_\E$ in \cite{MBays} is more precisely determined which is important in the context of this paper.) 

We may assume that $\Sigma_\A$ is stronger than  the first order theory
$\mathrm{Th}_\A$ in the same language with parameters $C=\{ c_a: a\in \kk\},$ and so
any Zariski closed subset definable over $\kk$ is first-order definable  in this theory. Recall the following folklore fact mentioned e.g. in \cite{Z1} and proved in detail in \cite{MBays}, Fact A.21. 

\medskip

{\bf Fact.} Let $\F$ be an algebraically closed field of characteristic 0 and $\X(\F)$ be $\F$-points an algebraic variety over $\kk\subset \F$ considered as the structure with relations corresponding to $\kk$-definable Zariski closed subsets of $\X^m(\F),$ $m\in \N.$ Then $(\F; +,\cdot, c_a)_{a\in \kk}$
is definable in the structure.

Consequently, if $\X'$ is a variety over $\kk'\subset \F$ and
 $\X'(\F)$ is a model of the theory $\mathrm{Th}_\X$ 
 then there is an automorphism $\sigma$ of field $\F$ such that $\kk^\sigma=\kk$ and $\X^\sigma(\F)=\X'(\F).$
  
\medskip

Let $\X_1$ and $\X_2$ be two complex algebraic  varieties over $\kk_1$ and $\kk_2$ respectively.
We will say that $\X_1$ and $\X_2$ have the {\bf same $L_{\omega_1,\omega}$-type} if  there is an uncountably categorical $L_{\omega_1,\omega}$-sentence $\Sigma$ over parameters $C$ interpreted as elements of $\kk_1$ in $\X_1$ and of $\kk_2$ in $\X_2$ such that the cover structures of both $\X_1(\C)$ and $\X_2(\C)$ are models of $\Sigma.$

\medskip

{\bf Theorem.} {\em  Let $\A_1$ and $\A_2$ be two complex abelian varieties.

Suppose $\A_1$ and $\A_2$  have the  same $L_{\omega_1,\omega}$-type.
Then there is a biholomorphic isomorphism $$\varphi: \A_1(\C)\to \A_2(\C),\mbox{ or }\varphi: \A_1(\C)\to \A_2^c(\C)$$
where $\A_2^c$ is the elliptic curve obtained by complex conjugation from $\A_2.$
} 

{\bf Proof.} We assume that $\A_1$ is over $\kk_1$ and $\A_2$ over $\kk_2,$ subfields of $\C.$
We will write $\U(\C)$ for the cover sort for both 
abelian varieties, and $\pi_1: \U\to \A_1,$ $\pi_2: \U\to \A_2$ the respective covers.

By categoricity and the Fact above there is an abstract isomorphism
$\phi$ between the cover structures, which have the restrictions to subsets, 
$$\phi_\kk:   \kk_1\to \kk_2;$$
  
 $$\phi_\Lambda: \la \Lambda_1\ra \to \la \Lambda_2\ra$$
where  $\la \Lambda_1\ra=\Lambda_1\otimes \Q,$  $\la \Lambda_2\ra= \Lambda_2\otimes \Q$
and 
$$\phi_{\pi\Lambda}:  \pi_1\la \Lambda_1\ra\to \pi_2\la \Lambda_2\ra.$$

 
For   $\kk_1$-interpretable varieties    $W$ and $V$ in $\A_1$ and a    
 $\kk_1$-rational map $f: W\to V$ we define $f^\phi: W^\phi\to V^\phi$ the $\kk_2$-rational map between the varieties  $W^\phi,$ $ V^\phi$ $\kk_2$-interpretable  in $\A_2,$ given by application of $\phi.$

Note that 
$\la \Lambda_1\ra$ is a topological metric space in the metric of $\Q$-vector spaces. It is 
metrically dense in $\U(\C)$ and so 
is $\pi_1\la \Lambda_1\ra$  in $\A_1(\C),$ and the same holds for $\A_2$ and its cover structure. 

Note also that $\phi_\Lambda$ is continuous and so is $\phi_{p\Lambda},$ since $\pi_1,$ $\pi_2$ are continuous.

Now we define $$\bar{\phi}:\left\lbrace \begin{array}{ll}
\ \U(\C) \to \U(\C)\\
\A_1(\C) \to \A_2(\C)
\end{array},\right.$$
the maps which are defined by continuous extension of $\phi_\Lambda$ and
of $\phi_{p\Lambda},$
that is, for $z=\lim z_i,$ $z_i\in \la \Lambda_1\ra,$ 
$$\bar{\phi}(z)= \lim \phi_\Lambda(z_i)$$
and for $z=\lim z_i,$ $z_i\in \pi_1\la \Lambda_1\ra,$ 
$$\bar{\phi}(z)= \lim \phi_{p\Lambda}(z_i).$$
Thus $\bar{\phi}$ is a continuous map between the two complex structures  the restrictions $\bar{\phi}_\Lambda$ and $\bar{\phi}_{\pi \Lambda}$ of which coinside with the restrictions of the isomorphism $\phi.$
We are going to prove that it is an isomorphism of the cover structures.

\medskip

Let $q: \A_1(\C)\to \C$ be a regular non-constant map over $\kk_1,$ defined on a Zariski open subset $D_{1,q}\subset \A_1(\C).$ 
 Let $q^\phi:\A_2(\C)\to \C $ be the respective map over $\kk_2.$

Note that $\pi_1\la \Lambda_1\ra\cap D_{1,q}$ is dense in $D_{1,q}$ and so
$q( \pi_1\la \Lambda_1\ra)$ is dense in $\C$ since $q$ is continuous and its image contains a  Zariski open  subset of $\C.$ 
Similarly,
 $q^\phi( \pi_2\la \Lambda_2\ra)$  is dense in $\C.$  Let $$\phi_q: q( \pi_1\la \Lambda_1\ra)\to q^\phi( \pi_2\la \Lambda_2\ra) $$ 
be the restriction of the isomorphism $\phi$ on the subsets.

The map $\phi_q$ is a homeomorphism since 
$\phi_{\pi\Lambda}$ is. 

 Thus the restriction of the isomorphism $\phi$ on the subsets of the form $q( \pi_1\la \Lambda_1\ra)$ is conituous.

  Define
$\bar{\phi}$ on $\C$  by extending $\phi_q$ continuously. This also defines $\bar{\phi}$ on $\C.$

If $q': \A_1\to \C$ is another such function then $\phi$ is continuous on $q( \pi_1\la \Lambda_1\ra)\cup q'( \pi_1\la \Lambda_1\ra).$
Define $\bar{\phi}'$ by extending $\phi$ from the union continuously. We will have  $\bar{\phi}'=\bar{\phi},$ since the two coinside on the dense subset $q( \pi_1\la \Lambda_1\ra).$
Thus, $\bar{\phi}$ does not depend on the choice of $q.$

\medskip

More generally, using the fact that $q( \pi_1\la \Lambda_1\ra)\subseteq \dcl(\Lambda_1)$ we get:
$$\bar{\phi}: \C\to \C$$
is the unique continuous map such that for any $\omega\in \dcl(\Lambda_1)\cap \C$
$$\bar{\phi}(\omega)=\phi(\omega).$$

\medskip

{\bf Claim.} Let $g: \C^m\to \C$ be a $\kk_1$-rational map. Then, for all $z\in \C^m,$
$$\bar{\phi}(g(z))=g^\phi (\bar{\phi}(z))).$$

Proof. Let $z=\lim z_i$ for some sequence $z_i\in \dcl(\Lambda_1).$ Since $g$ is continuous, $g(z)=\lim g(z_i)$ and since 
$\bar{\phi}$ is continuous on $\C$
$$\bar{\phi}(g(z))=\lim \bar{\phi}(g(z_i))=\lim \phi(g(z_i))=\lim g^\phi(\phi(z_i))=g^\phi(\lim \phi(z_i))=
g^\phi (\bar{\phi}(z)))$$

\medskip

{\bf Corollary}. $\bar{\phi}$ is a continuous isomorphism of $\C,$ that is an identity or complex conjugation $z\mapsto z^c.$ In particular
$\phi_\kk: \kk_1\to \kk_2$ is  an identity or complex conjugation.

Indeed, take $g(z_1,z_2):=z_1*z_2$ where $*$ is + or $\cdot$. 
 Then $$\bar{\phi}(z_1* z_2)=\bar{\phi}(z_1)*\bar{\phi}(z_2).$$

 $\Box$

\section{Shimura varieties}

\medskip

 Let $\X_1,$ $\X_2$ be  connected Shimura varieties (see \cite{Milne0}, Def.4.4) over their reflex fields $\kk_1,\kk_2$ given by the Shimura datum $(\G_1,\mathcal{H}_1),$ $(\G_2,\mathcal{H}_2)$ and   congruence subgroups $\Gamma_1\subset \G_1(\Q),$  $\Gamma_2\subset \G_2(\Q),$ respectively.
$$\X_1=\Gamma_1\backslash \mathcal{H}_1\mbox{ and }\X_2=\Gamma_2\backslash \mathcal{H}_2$$
with respective covering maps
$$j_1: \mathcal{H}_1\to \X_1\mbox{ and }j_2: \mathcal{H}_2\to \X_2.$$
 
 
A Shimura datum $(\G,\mathcal{H}),$ also associates with each point $h\in \mathcal{H}_1$ 
a homomorphisms $h: U(1) \to \G^\mathrm{ad}_\R$ from the $1$-dimensional real torus $U(1),$ and all the homomorphisms are linked by conjugation by elements of
  $\G^\mathrm{ad}(\R)^+.$ The image $h(U(1))=T_h$ is a torus subgroup of 
  $\G^\mathrm{ad}_\R$ which can be identified as a subgroup fixing the point $h.$ Since there is only one non-identity continuous automorphism of $U(1),$ $u\mapsto u\inv,$
  there are exactly two $\G^\mathrm{ad}_\R$-conjugacy classes of homomorphisms $h: U(1) \to \G^\mathrm{ad}_\R$ which can be called
  $\mathcal{H}^+$ and  $\mathcal{H}^-.$ 
  
  It follows from Shimura theory, see \cite{BG}, that
   $\X,$  a complex Shimura variety and $\X^c$ its complex conjugate, 
  
 \be\label{c} \X=\Gamma\backslash \mathcal{H}^+ \Leftrightarrow \X^c=\Gamma\backslash \mathcal{H}^-\ee
 and the action of complex conjugation corresponds to the automorphism $u\to u\inv$ on $U(1).$
 
 \medskip
 
 Below we use the notation $\mathcal{H}_1$ and $\mathcal{H}_2$ for sets and  $\mathcal{H}_1^{\pm}$ and $\mathcal{H}_2^{\pm}$ for their presentations as classes of homomomorphisms from $U(1).$
 
 \medskip
 
 {\bf Fact} (Th 5.4,\cite{Milne0}). Let $\G$ be a connected algebraic group over $\Q.$
Then $\G(\Q)$ is dense in $\G(\R).$ 

\medskip
 
{\bf Theorem.} {\em Suppose there is a categorical  $L_{\omega_1,\omega}$-sentence $\Sigma$ over parameters $C$ such that both structures $(\G_1,\mathcal{H}_1, \Gamma_1, \X_1, j_1)$ and
 $(\G_2,\mathcal{H}_2, \Gamma_2, \X_2,j_2)$ are models of $\Sigma$ when one interprets $C$ as $\kk_1$ and $\kk_2$ respectively. Then 
there is a biholomorphic isomorphism} \be\label{var} \varphi: \X_1(\C)\to \X_2(\C),\mbox{ or }\varphi: \X_1(\C)\to \X_2^c(\C).\ee

{\bf Proof.} First we note that 
it follows from the categoricity assumption and the Fact of section \ref{s1} that there exists an abstract isomorphism $\phi$ between the two structures which is realised as ismorphisms

\be \label{k} \phi:\kk_1\to \kk_2\ee
\be \label{X} \phi: \X_1\to \X_2\ee
\be \label{H} \phi: \mathcal{H}_1\to \mathcal{H}_2\ee
\be \label{GQ}\phi: \G_1(\Q)\to \G_2(\Q);\ee 
 \be \label{Gamma}\phi: \Gamma_1\to \Gamma_2\ee
\be\label{gx} \forall g\in \G_1(\Q), \forall x\in  \mathcal{H}_1\ \
\phi (g\cdot x)=\phi(g)\cdot \phi(x) \ee

and 
\be \label{j} j_2\circ \phi=\phi\circ j_1.\ee

Since $\G_1(\Q)$ is dense in $\G_1(\R)$ and  $\G_2(\Q)$  in $\G_2(\R)$
there is a unique extension of $\phi$ to the continuous isomorphism 
$$\phi_\R: \G_1(\R)\to \G_2(\R), \ \G^\mathrm{ad}1(\R)\to \G^\mathrm{ad}2(\R)$$

together with the isomorphism of the action (\ref{gx}) since it is continuous both in $g$ and $x:$
$$ \forall g\in \G_1(\R), \forall x\in  \mathcal{H}_1\ \
\phi_\R (g\cdot x)=\phi_\R(g)\cdot \phi_\R(x)$$
(where $\G^\mathrm{ad}(\R)=\G(\R)/Z,$ $Z$ subgroup acting trivially).


Thus, $\phi$ induces 
$$\phi_{\pm}: \{ h: U(1)\to \G_1^\mathrm{ad}(\R) \} \to 
\{ \phi_{\pm}(h): U(1)\to \G_2^\mathrm{ad}(\R)\}$$

Since $h\inv\circ\phi_\R\inv\circ\phi_{\pm}(h): U(1)\to U(1)$ is either identity or $u\to u\inv,$ there are two possibilities:

$$\phi_{\pm}:  \mathcal{H}_1^+\to \mathcal{H}_2^+$$
and 

$$\phi_{\pm}:  \mathcal{H}_1^-\to \mathcal{H}_2^+$$

 Respectively, by fact (\ref{c}), we have two biholomorphic isomorphisms  (\ref{var}). $\Box$

\end{document}